\documentclass[aps,a4paper]{revtex4}
\usepackage{epsfig}
\usepackage{graphicx}
\usepackage{amsmath,amssymb,color}
\usepackage[english]{babel}

\parskip=\medskipamount

\newcommand{\eq}[1]{(\ref{#1})}
\newcommand{\fig}[1]{Fig.~\ref{#1}}

\newcommand{\be}{\begin{equation}}
\newcommand{\ee}{\end{equation}}

\newcommand\disp{\displaystyle}

\newcommand{\la}{\left<}
\newcommand{\ra}{\right>}
\newcommand{\eps}{\varepsilon}

\newcommand{\re}{\textrm{Re}\,}
\newcommand{\im}{\textrm{Im}\,}

\begin{document}

\title{Number-theoretic aspects of 1D localization: "popcorn function" with Lifshitz tails
and its continuous approximation by the Dedekind eta}

\author{S. Nechaev$^{1,2}$, and K. Polovnikov$^{3,4}$}

\address{$^1$Interdisciplinary Scientific Center Poncelet (ISCP), Bolshoy Vlasyevskiy Pereulok 11,
119002, Moscow, Russia \\ $^2$P.N. Lebedev Physical Institute RAS, 119991, Moscow, Russia \\
$^3$Center for Energy Systems, Skolkovo Institute of Science and Technology, 143005 Skolkovo,
Russia \\ $^4$Physics Department, M.V. Lomonosov Moscow State University, 119992 Moscow, Russia}

\begin{abstract}

We discuss the number-theoretic properties of distributions appearing in physical systems when an
observable is a quotient of two independent exponentially weighted integers. The spectral density
of ensemble of linear polymer chains distributed with the law $\sim f^L$ ($0<f<1$), where $L$ is
the chain length, serves as a particular example. At $f\to 1$, the spectral density can be
expressed through the discontinuous at all rational points, Thomae ("popcorn") function. We suggest
a continuous approximation of the popcorn function, based on the Dedekind $\eta$-function near the
real axis. Moreover, we provide simple arguments, based on the "Euclid orchard" construction, that
demonstrate the presence of Lifshitz tails, typical for the 1D Anderson localization, at the
spectral edges. We emphasize that the ultrametric structure of the spectral density is ultimately
connected with number-theoretic relations on asymptotic modular functions. We also pay attention to
connection of the Dedekind $\eta$-function near the real axis to invariant measures of some
continued fractions studied by Borwein and Borwein in 1993.

\end{abstract}

\maketitle

\section{Introduction}

The so-called "popcorn function" \cite{popcorn}, $g(x)$, known also as the Thomae function, has
also many other names: the raindrop function, the countable cloud function, the modified Dirichlet
function, the ruler function, etc. It is one of the simplest number-theoretic functions possessing
nontrivial fractal structure (another famous example is the everywhere continuous but never
differentiable Weierstrass function). The popcorn function is defined on the open interval $x \in
(0, 1)$ according to the following rule:
\be
g(x) = \begin{cases} \frac{1}{q} & \mbox{if $x=\frac{p}{q}$, and $(p,q)$ are coprime} \medskip \\
0 & \mbox{if $x$ is irrational} \end{cases}
\ee
The popcorn function $g$ is discontinuous at every rational point because irrationals come
infinitely close to any rational number, while $g$ vanishes at all irrationals. At the same time,
$g$ is continuous at irrationals.

\begin{figure}[ht]
\centerline{\includegraphics[width=10cm]{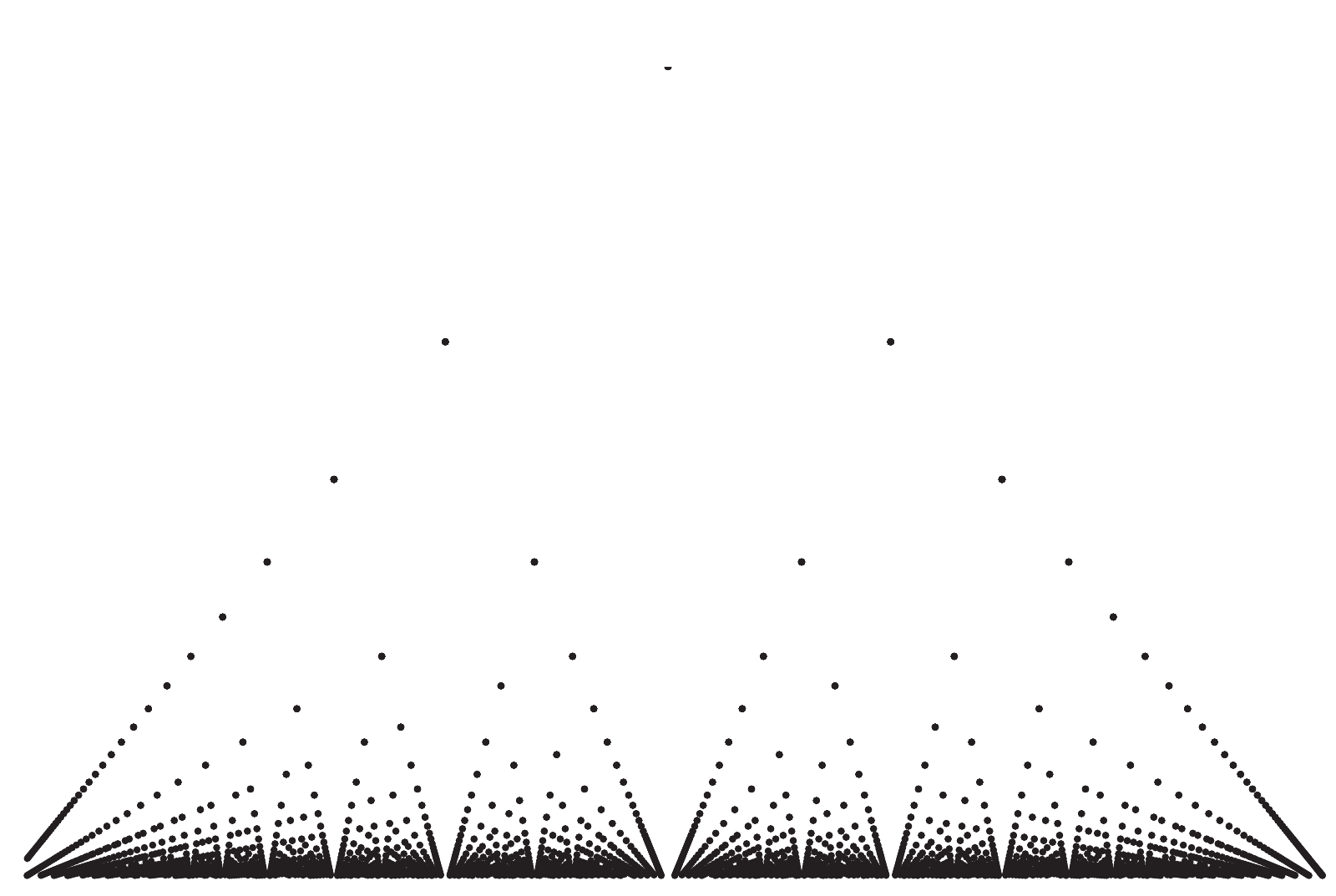}}
\caption{Popcorn (Thomae) function.}
\label{fig:01}
\end{figure}

In order to demonstrate this fact accurately, consider some irrational number, $t$, at which
$g(t)=0$, and take some $0 < \eps < 1$. Without the loss of generality $\eps$ is assumed to be
rational, otherwise one may replace $\eps$ with any smaller rational $\eps'=\frac{k}{m}<\eps$ such
that $\mathrm{gcd}(k, m) = 1$. Thus, there is a finite set of the rational numbers $\Omega_\eps =
\{q= \frac{i}{j}, 1<j \le m, 1 \le i < j\}$, whose popcorn-values are not smaller than $g(\eps)$.
Now assign $\delta(\eps) = inf\{|q - t|, q \in \Omega_\eps\}$ that defines the vicinity of $t$, in
which the values of $g$ are smaller than $\eps$:
\be
|t-y|<\delta(\eps) \rightarrow |g(t)-g(y)|=|g(y)|<\eps
\ee
These inequalities prove the continuity of of the popcorn function $g(x)$ at irrational values of
$x$.

To our point of view, the popcorn function has not yet received decent attention among researchers,
though its emergence in various physical problems seems impressive, as we demonstrate below.
Apparently, the main difficulty deals with the discontinuity of $g(x)$ at every rational point,
which often results in a problematic theoretical treatment and interpretation of results for the
underlying physical system. Thus, a natural, physically justified "continuous approximation" to
the popcorn function is very demanded.

Here we provide such an approximation, showing the generality of the "popcorn-like" distributions
for a class of one-dimensional disordered systems. From the mathematical point of view, we
demonstrate that the popcorn function can be constructed on the basis of the modular Dedekind
function, $\eta(x+iy)$, when the imaginary part, $y$, of the modular parameter $z=x+iy$ tends to 0.

One of the most beautiful incarnations of the popcorn function arises in a so-called "Euclid
orchard" representation. Consider an orchard of trees of \emph{unit hights} located at every point
$(an, am)$ of the two-dimensional square lattice, where $n$ and $m$ are nonnegative integers
defining the lattice, and $a$ is the lattice spacing, $a = 1/\sqrt{2}$. Suppose we stay on the line
$n = 1 - m$ between the points $A(0, a)$ and $B(a, 0)$, and observe the orchard grown in the first
quadrant along the rays emitted from the origin $(0,0)$ -- see the \fig{fig:01}.

\begin{figure}[ht]
\centerline{\includegraphics[width=15cm]{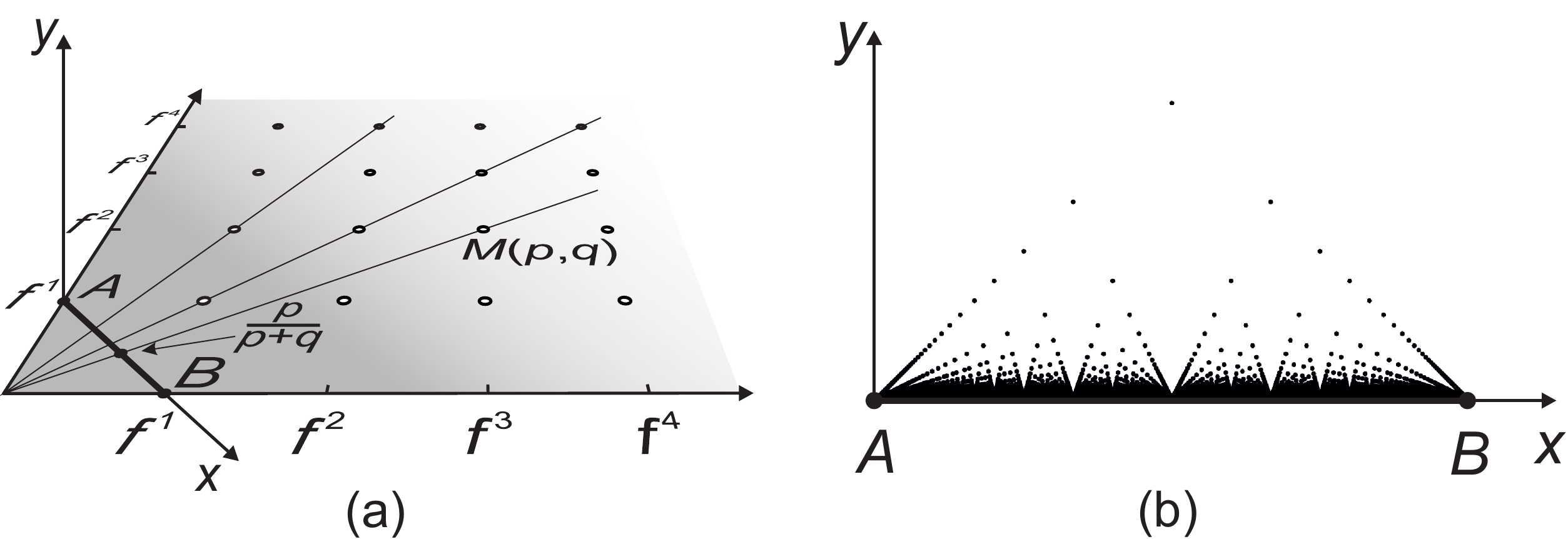}}
\caption{(a) Construction of the Euclid Orchard; (b) Thomae function.}
\label{fig:01}
\end{figure}

Along these rays we see only the first open tree with coprime coordinates, $M(ap,aq)$, while all
other trees are shadowed. Introduce the auxiliary coordinate basis $(x,y)$ with the axis $x$ along
the segment $AB$ and $y$ normal to the orchard's plane (as shown in the \fig{fig:01}a). We set the
origin of the $x$ axis at the point $A$, then the point $B$ has the coordinate $x = a$. It is a
nice school geometric problem to establish that: (i) having the focus located at the origin, the
tree at the point $M(ap,aq)$ is spotted at the place $x = \frac{p}{p+q}$, (ii) the \emph{visible}
height of this tree is $\frac{1}{p+q}$. In other words, the "visibility diagram" of such a lattice
orchard is exactly the popcorn function.

The popcorn correspondence $\frac{p}{q} \to \frac{1}{q}$ arises in the Euclid orchard problem as a
purely geometrical result. However, the same function has appeared as a probability distribution in
a plethora of biophysical and fundamental problems, such as the distribution of quotients of reads
in DNA sequencing experiment \cite{dna}, quantum $1/f$ noise and Frenel-Landau shift \cite{planat},
interactions of non-relativistic ideal anyons with rational statistics parameter in the magnetic
gauge approach \cite{lundholm}, or frequency of specific subgraphs counting in the protein-protein
network of a \textit{Drosophilla} \cite{drosophilla}. Though the extent of similarity with the
original popcorn function could vary, and experimental profiles may drastically depend on
peculiarities of each particular physical system, a general probabilistic scheme resulting in the
popcorn-type manifestation of number-theoretic behavior in nature, definitely survives.

Suppose two random integers, $\phi$ and $\psi$, are taken independently from a discrete probability
distribution, $Q_n = f^n$, where $f = 1 - \eps > 0$ is a "damping factor". If $\mathrm{gcd}(p,q) =
1$, then the combination $\nu  = \frac{\phi}{\phi+\psi}$ has the popcorn-like distribution $P(\nu)$
in the asymptotic limit $\eps \ll 1$:
\be
P\left(\nu = \frac{p}{p+q}\right) = \sum_{n=1}^{\infty} f^{n(p+q)}
= \frac{(1-\eps)^{p+q}}{1 - (1-\eps)^{p+q}} \approx \frac{1}{\eps(p+q)}
\label{ind}
\ee

The formal scheme above can be understood on the basis of the generalized Euclid orchard
construction, if one considers a $(1+1)$-dimensional directed walker on the lattice (see
\fig{fig:01}a), who starts from the lattice origin and makes $\phi$ steps along one axis, followed
by $\psi$ steps along the other axis. At every step the walker can die with the probability $\eps =
1 - f$. Then, having an ensemble of such walkers, one arrives at the "orchard of walkers", i.e. for
some point $\nu$ on the $x$ axis, a fraction of survived walkers, $P(\nu)$, is given exactly by the
popcorn function.

In order to have a relevant physical picture, consider a toy model of diblock-copolymer
polymerization. Without sticking to any specific polymerization mechanism, consider an ensemble of
diblock-copolymers $AB$, polymerized independently from both ends in a cloud of monomers of
relevant kind (we assume, only $A-A$ and $B-B$ links to be formed). Termination of polymerization
is provided by specific "radicals" of very small concentration, $\eps$: when a radical is attached
to the end (irrespectively, $A$ or $B$), it terminates the polymerization at this extremity
forever. Given the environment of infinite capacity, one assigns the probability $f = 1-\eps$ to a
monomer attachment at every elementary act of the polymerization. If $N_A$ and $N_B$ are molecular
weights of the blocks $A$ and $B$, then the composition probability distribution in our ensemble,
$P\left(\varphi = \frac{N_A}{N_A+N_B}\right)$, in the limit of small $\eps \ll 1$ is "ultrametric"
(see \cite{avetisov} for the definition of the ultrametricity) and is given by the popcorn
function:
\be
P\left(\varphi = \frac{p}{p + q}\right) \approx \frac{1}{\eps(p+q)}
\stackrel{def} = \frac{1}{\eps}g(\varphi)
\ee

In the polymerization process described above we have assumed identical independent probabilities
for the monomers of sorts ("colors") $A$ and $B$ to be attached at both chain ends. Since no
preference is implied, one may look at this process as at a homopolymer ("colorless") growth,
taking place at two extremities. For this process we are interested in statistical characteristics
of the resulting ensemble of the homopolymer chains. What would play the role of "composition" in
this case, or in other words, how should one understand the fraction of monomers attached at one
end? As we show below, the answer is rather intriguing: the respective analogue of the probability
distribution is the spectral density of the ensemble of linear chains with the probability $Q_L$
for the molecular mass distribution, where $L$ is the length of a chain in the ensemble.

\section{Spectral statistics of exponentially weighted ensemble of linear graphs}

\subsection{Spectral density and the popcorn function}

The former exercises are deeply related to the spectral statistics of ensembles of linear polymers.
In a practical setting, consider an ensemble of noninteracting linear chains with exponential
distribution in their lengths.  We claim the emergence of the fractal popcorn-like structure in the
spectral density of corresponding adjacency matrices describing the connectivity of elementary
units (monomers) in linear chains.

The ensemble of exponentially weighted homogeneous chains, is described by the bi-diagonal
symmetric $N\times N$ adjacent matrix $B=\{b_{ij}\}$:
\be
B = \left(\begin{array}{ccccc}
0 & x_1 & 0 & 0 & \cdots \\  x_1 & 0 & x_2 & 0 &  \\  0 & x_2 & 0 & x_3 &  \\
0 & 0 & x_3 & 0 &  \\ \vdots &  &  &  & \ddots
\end{array} \right)
\label{eq:06}
\ee
where the distribution of each $b_{i,i+1}=b_{i+1,i}=x_i$ ($i=1,...,N$) is Bernoullian:
\be
x_i=\left\{\begin{array}{ll} 1 & \mbox{with probability $x$} \medskip \\
0 & \mbox{with probability $\eps = 1-f$} \end{array} \right.
\label{eq:06a}
\ee
We are interested in the spectral density, $\rho_{\eps}(\lambda)$, of the ensemble of matrices $B$
in the limit $N\to\infty$. Note that at any $x_k=0$, the matrix $B$ splits into independent blocks.
Every $n\times n$ block is a symmetric $n\times n$ bi-diagonal matrix $A_n$ with all $x_k=1$,
$k=1,...,n$, which corresponds to a chain of length $n$. The spectrum of the matrix $A_n$ is
\be
\lambda_{k,n} = 2\cos\frac{\pi k}{n+1}; \qquad (k=1,...,n)
\label{eq:06b}
\ee

All the eigenvalues $\lambda_{k,n}$ for $k=1,...,n-1$ appear with the probability $Q_n = f^n$ in
the spectrum of the matrix \eq{eq:06}. In the asymptotic limit $\eps \ll 1$, one may deduce an
equivalence between the composition distribution in the polymerization problem, discussed in the
previous section, and the spectral density of the linear chain ensemble. Namely, the probability of
a composition $\varphi = \frac{p}{p+q}$ in the ensemble of the diblock-copolymers can be precisely
mapped onto the peak intensity (the degeneracy) of the eigenvalue $\lambda=\lambda_{p,p+q-1} =
2\cos\frac{\pi p}{p+q}$ in the spectrum of the matrix $B$. In other words, the integer number $k$
in the mode $\lambda_{k,n}$ matches the number of $A$-monomers, $N_A = kz$, while the number of
$B$-monomers matches $N_B = (n + 1 - k)z$, where $z \in N$, in the respective diblock-copolymer.

The spectral statistics survives if one replaces the ensemble of Bernoullian two-diagonal adjacency
matrices $B$ defined by \eq{eq:06}--\eq{eq:06a} by the ensemble of random Laplacian matrices.
Recall that the Laplacian matrix, $L=\{a_{ij}\}$, can be constructed from
adjacency matrix, $B=\{b_{ij}\}$, as follows: $a_{ij} = -b_{ij}$ for $i\neq j$, and $a_{ii} =
\sum_{j=1}^{N}b_{ij}$. A search for eigenvalues of the Laplacian matrix $L$ for linear chain, is
equivalent to determination its relaxation spectrum. Thus, the density of the relaxation spectrum
of the ensemble of noninteracting linear chains with the exponential distribution in lengths, has
the signature of the popcorn function.

\begin{figure}[ht]
\centerline{\includegraphics[width=8cm]{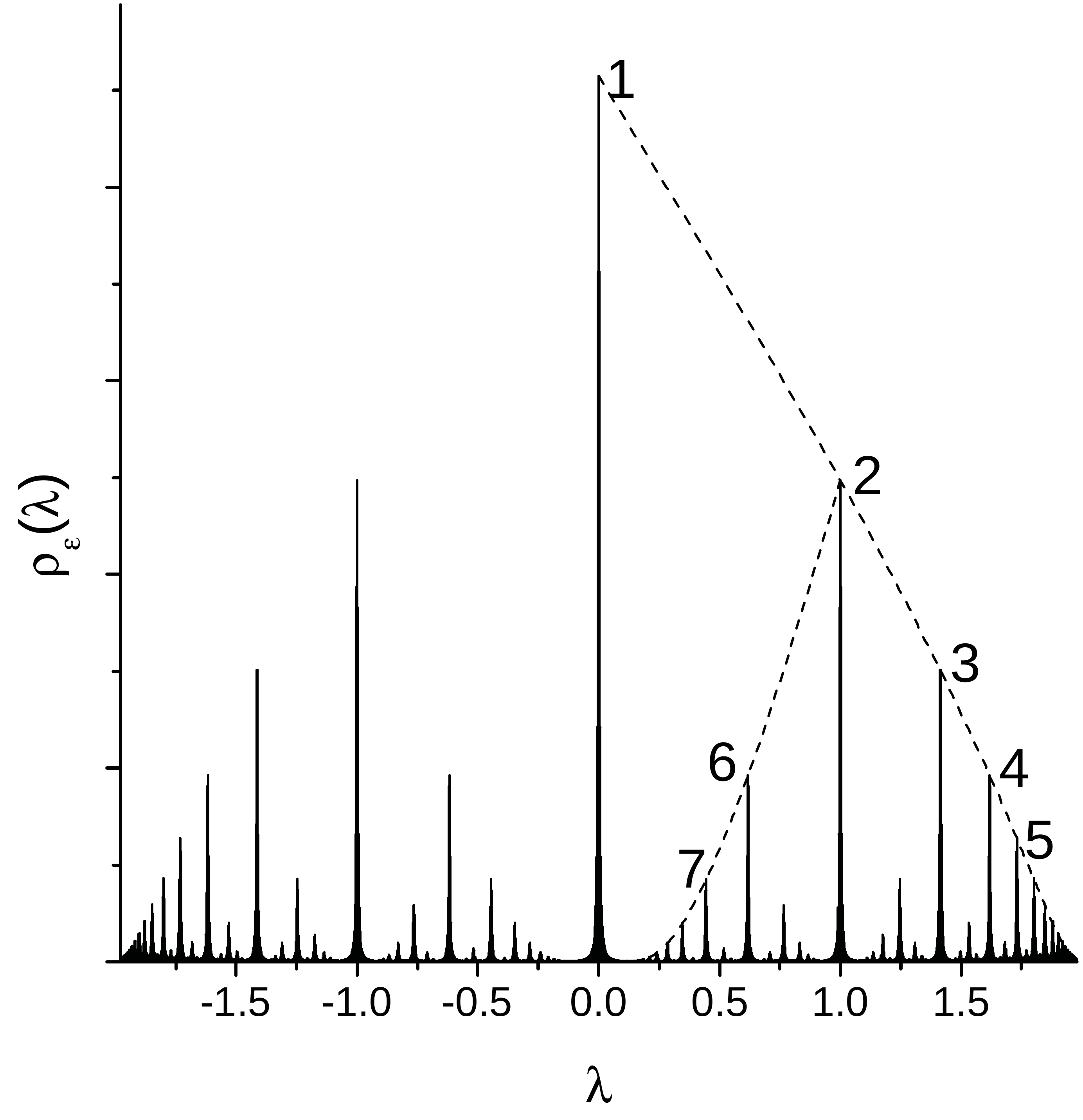}}
\caption{The spectral density $\rho_{\eps}(\lambda)$ for the ensemble of bi-diagonal matrices of
size $N=10^3$ at $f=0.7$. The regularization parameter $\eps$ is taken $\eps=2\times 10^{-3}$.}
\label{fig:02}
\end{figure}

To derive $\rho_{\eps}(\lambda)$ for arbitrary values of $\eps$, let us write down the spectral
density of the ensemble of $N\times N$ random matrices $B$ with the bimodal distribution of the
elements as a resolvent:
\be
\rho_\eps(\lambda) = \lim_{N\to\infty}\Big< \sum_{k=1}^{n} \delta(\lambda-\lambda_{kn}) \Big>_{Q_n}
= \lim_{N\to\infty \atop y\to+0} y\; \im\, \Big< G_n(\lambda - iy) \Big>_{Q_n}  = \lim_{N\to\infty
\atop y\to+0} y \sum_{n=1}^N Q_n\; \im\, G_n(\lambda - iy)
\label{rh}
\ee
where $\la ...\ra_{Q_n}$ means averaging over the distribution $Q_n=(1-\eps)^n$, and the following
regularization of the Kronecker $\delta$-function is used:
\be
\delta(\xi) = \lim_{y\to+0} \im \frac{y}{\xi- iy}
\label{delt}
\ee
The function $G_n$ is associated with each particular gapless matrix $B$ of $n$ sequential "1" on
the sub-diagonals,
\be
G_n(\lambda-iy) = \sum_{k=1}^n \frac{1}{\lambda-\lambda_{k,n}-iy}
\label{green}
\ee
Collecting \eq{eq:06b}, \eq{rh} and \eq{green}, we find an explicit expression for the density of
eigenvalues:
\be
\rho_\eps(\lambda) = \lim_{N\to\infty \atop y\to+0} y \sum_{n=1}^{N} (1-\eps)^n
\sum_{k=1}^n\frac{y}{\left(\lambda-2\cos\frac{\pi k}{n+1}\right)^2+y^2}
\label{eq:10}
\ee
The behavior of the inner sum in the spectral density in the asymptotic limit $y \to 0$ is easy to
understand: it is $\frac{1}{y}$ at $\lambda=2\cos{\frac{\pi k}{n+1}}$ and zero otherwise. Thus, one
can already infer a qualitative similarity with the popcorn function. It turns out, that the
correspondence is quantitative for $\eps = 1 - f \ll 1$. Driven by the purpose to show it, we
calculate the values of $\rho_\eps(\lambda)$ at the peaks, i.e. at rational points $\lambda =
2\cos{\frac{\pi p}{p+q}}$ with $\mathrm{gcd}(p, q) = 1$ and end up with the similar geometrical
progression, as for the case of diblock-copolymers problem \eq{ind}:
\be
\rho_\eps \left(\lambda = 2\cos{\frac{\pi p}{p+q}}\right) = \sum_{s=1}^{\infty} (1-\eps)^{(p+q)s-1}
= \frac{(1-\eps)^{p+q-1}}{1 - (1-\eps)^{p+q}}\Bigg|_{\eps\to 0} \approx \frac{1}{\eps(p+q)}
\stackrel{def} = g\left(\frac{1}{\pi}\arccos{\frac{\lambda}{2}}\right)
\label{eq:115}
\ee
The typical sample plot $\rho_{\eps}(\lambda)$ for $f=0.7$ computed numerically via \eq{eq:10}
with $\eps=2\times 10^{-3}$ is shown in the \fig{fig:02} for $N=10^3$.

\subsection{Enveloping curves and tails of the eigenvalues density}

\begin{figure}[ht]
\centerline{\includegraphics[width=15cm]{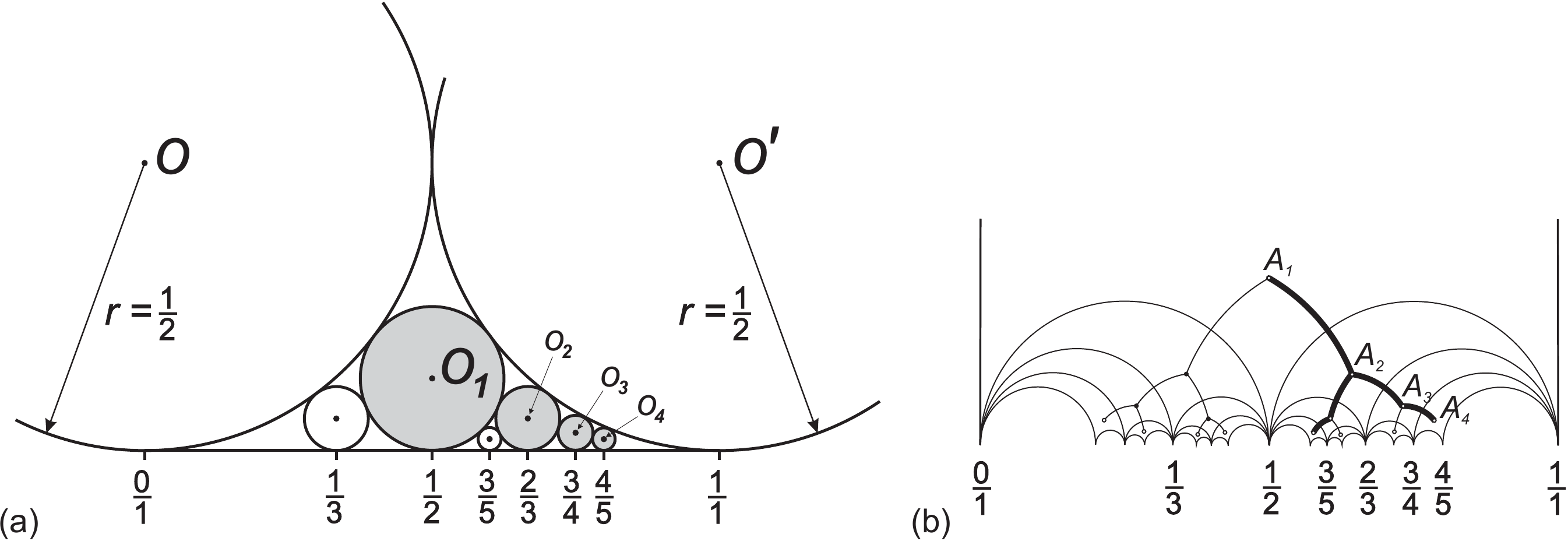}}
\caption{Ford circles as illustration of the Farey sequence construction: (a) Each circle touches
two neighbors (right and left) and the segment. The position of newly generated circle is
determined via the $\oplus$ addition: $\frac{p_{k-1}}{q_{k-1}} \oplus \frac{p_{k+1}}{q_{k+1}}=
\frac{p_{k-1}+p_{k+1}}{q_{k-1}+q_{k+1}}$; (b) The same Farey sequence generated by sequential
fractional-linear transformations of the fundamental domain of the modular group $SL(2,Z)$. }
\label{fig:farey}
\end{figure}

Below we pay attention to some number-theoretic properties of the spectral density of the argument
$-\lambda$, since in this case the correspondence with the composition ratio is precise. One can
compute the enveloping curves for any monotonic sequence of peaks depicted in \fig{fig:02}, where
we show two series of sequential peaks: $S_1=\{$1--2--3--4--5--...$\}$ and
$S_2=\{$2--6--7--...$\}$. Any monotonic sequence of peaks corresponds to the set of eigenvalues
$\lambda_{k,n}$ constructed on the basis of a Farey sequence \cite{farey}. For example, as shown
below, the peaks in the series $S_1$ are located at:
$$
\lambda_k = -\lambda_{k,k} =-2\cos\frac{\pi k}{k+1}, \qquad (k=1,2,...)
$$
while the peaks in the series $S_2$ are located at:
$$
\lambda_{k'} = -\lambda_{k',2k'-2} = -2\cos \frac{\pi k'}{2k'-1}, \qquad (k'=2,3,...)
$$
Positions of peaks obey the following rule: let $\{\lambda_{k-1},\, \lambda_k,\, \lambda_{k+1}\}$
be three consecutive monotonically ordered peaks (e.g., peaks 2--3--4 in \fig{fig:02}), and let
$$
\lambda_{k-1}=-2\cos \frac{\pi p_{k-1}}{q_{k-1}}, \quad \lambda_{k+1}=-2\cos \frac{\pi
p_{k+1}}{q_{k+1}}
$$
where $p_k$ and $q_k$ ($k=1,...,N$) are coprimes. The position of the intermediate peak,
$\lambda_k$, is defined as
\be
\lambda_{k}=-2\cos \frac{\pi p_{k}}{q_{k}}; \qquad \frac{p_{k}}{q_{k}} = \frac{p_{k-1}}{q_{k-1}}
\oplus \frac{p_{k+1}}{q_{k+1}} \equiv \frac{p_{k-1}+p_{k+1}}{q_{k-1}+q_{k+1}}
\label{eq:11}
\ee
The sequences of coprime fractions constructed via the $\oplus$ addition are known as Farey
sequences. A simple geometric model behind the Farey sequence, known as Ford circles \cite{ford,
coxeter}, is shown in \fig{fig:farey}a. In brief, the construction goes as follows. Take the
segment $[0,1]$ and draw two circles $O_1$ and $O_2$ both of radius $r=\frac{1}{2}$, which touch
each other, and the segment at the points 0 and 1. Now inscribe a new circle $O_3$ touching $O_1$,
$O_2$ and $[0,1]$. Where is the position of the new circle along the segment? The generic recursive
algorithm constitutes the Farey sequence construction. Note that the same Farey sequence can be
sequentially generated by fractional-linear transformations (reflections with respect to the arcs)
of the fundamental domain of the modular group $SL(2,Z)$ -- the triangle lying in the upper
halfplane $\im z>0$ of the complex plane $z$ (see \fig{fig:farey}b).

Consider the main peaks series, $S_1=\{$1--2--3--4--5--...$\}$. The explicit expression for their
positions reads as:
\be
\lambda_k = -2\cos\frac{\pi k}{k+1}; \qquad k = 1, 2, ...
\ee

One can straightforwardly investigate the asymptotic behavior of the popcorn function in the limit
$k \to \infty$. From \eq{eq:115} one has for arbitrary $f < 1$ the set of parametric equations:
\be
\left\{\begin{array}{l}
\disp \rho_\eps(\lambda_k) = \frac{f^{k}}{1 - f^{k+1}}\Big|_{k\gg 1} \approx f^{k} \\
\disp \lambda_k=-2\cos\frac{\pi k}{k+1}\Big|_{k\gg 1} \approx 2 - \frac{\pi^2}{k^2}
\end{array}\right.
\label{asm}
\ee
From the second equation of \eq{asm}, we get $k\approx \frac{\pi}{\sqrt{2-\lambda}}$. Substituting
this expression into the first one of \eq{asm}, we end up with the following asymptotic behavior of
the spectral density near the spectral edge $\lambda\to 2^-$:
\be
\rho_{\eps}(\lambda) \approx \exp\left(\frac{\pi \ln f}{\sqrt{2-\lambda}}\right) \qquad (0<f<1)
\label{asm2}
\ee
The behavior \eq{asm} is the signature of the Lifshitz tail typical for the 1D Anderson
localization:
\be
\rho_{\eps}(E)\approx e^{-C E^{-D/2}};
\ee
where $E= 2-\lambda$ and $D=1$.

\section {From popkorn to Dedekind $\eta$-function}

\subsection {Some facts about Dedekind $\eta$-function and related series}

The popcorn function has discontinuous maxima at rational points and continuous valleys at
irrationals. We show in this section, that the popcorn function can be regularized on the basis of
the everywhere continuous Dedekind function $\eta(x+iy)$ in the asymptotic limit $y\to 0$.

The famous Dedekind $\eta$-function is defined as follows:
\be
\eta(z)=e^{\pi i z/12}\prod_{n=0}^{\infty}(1-e^{2\pi i n z})
\label{eq:dedeta}
\ee
The argument $z=x+iy$ is called the modular parameter and $\eta(z)$ is defined for $\im z>0$ only.
The Dedekind $\eta$-function is invariant with respect to the action of the modular group
$SL(2,\mathbb{Z})$:
\be
\begin{array}{l}
\eta (z+1)=e^{\pi i z/12}\;\eta(z) \medskip \\ \eta\left(-\frac{1}{z}\right) = \sqrt{-i}\; \eta(z)
\end{array}
\label{2}
\ee
And, in general,
\be
\eta\left(\frac{az+b}{cz+d}\right) = \omega(a,b,c,d)\; \sqrt{c z + d}\; \eta(z)
\label{3}
\ee
where $ad-bc=1$ and $\omega(a,b,c,d)$ is some root of 24th degree of unity \cite{dedekind}.

It is convenient to introduce the following "normalized" function
\be
h(z) = |\eta(z)| (\im z)^{1/4}
\label{eq:15}
\ee
The analytic structure of $h(z)$ has been discussed in \cite{avetisov} in the context of the
\emph{ultrametric} landscape construction. In particular, the function $h(z)$ satisfies the duality
relation, which follows from \eq{2}--\eq{3}:
\be
h\left(\left\{\frac{m}{k}\right\} + iy\right) = h\left(\left\{\frac{n}{k}\right\}+\frac{i}{k^2y}
\right)
\label{eq:15a}
\ee
where $m n - k r = 1,\; \{k,m,n,r\}\in \mathbb{Z},\; (y>0)$ and $\left\{\frac{m}{k}\right\},
\left\{\frac{s}{k}\right\}$ denote fractional parts of corresponding quotients. In Appendix we use
Eq.\eq{eq:15a} to compute the values of $h(x+iy)$ at rational points $x$ near the real axis, i.e.
at $y\to 0$.

The real analytic Eisenstein series $E(z, s)$ is defined in the upper half-plane, $H = \{z: \im(z)
> 0\}$ for $\re(s) > 1$ as follows:
\be
E(z,s) = \frac{1}{2}\sum_{\{m,n\} \in \mathbb{Z}^2 \backslash \{0, 0\}} \frac{y^s}{|mz + n|^{2s}};
\qquad z = x+iy
\ee
This function can be analytically continued to all $s$-plane with one simple pole at $s=1$. Notably
it shares the same invariance properties on $z$ as the Dedekind $\eta$-function. Moreover, $E(s,
z)$, as function of $z$, is the $SL(2,\mathbb{Z})$--automorphic solution of the hyperbolic Laplace
equation:
$$
-y^2 \left(\frac{\partial^2}{\partial x^2}+\frac{\partial^2}{\partial y^2}\right)
E(z, s) = s(1-s)\; E(z, s)
$$

The Eisenstein series is closely related to the Epstein $\zeta$-function, $\zeta(s, Q)$, namely:
\be
\zeta(s, Q) = \sum_{\{m,n\} \in \mathbb{Z}^2 \backslash \{0, 0\}} \frac{1}{Q(m, n)^s} =
\frac{2}{d^{s/2}}E(z, s),
\label{epstein}
\ee
where $Q(m, n) = am^2 + 2bmn + cn^2$ is a positive definite quadratic form, $d=ac-b^2 > 0$, and
$\disp z = \frac{-b+i\sqrt d}{a}$. Eventually, the logarithm of the Dedekind $\eta$-function is
known to enter in the Laurent expansion of the Epstein $\zeta$-function. Its residue at $s=1$ has
been calculated by Dirichlet and is known as the first Kronecker limit formula
\cite{epstein,siegel,motohashi}. Explicitly, it reads at $s\to 1$:
\be
\zeta(s, Q) = \frac{\pi}{\sqrt d} \frac{1}{s-1} \\ + \frac{2\pi}{\sqrt d}\left(\gamma +
\ln\sqrt{\frac{a}{4d}} - 2\ln|\eta(z)|\right) + O(s-1)
\label{eq:zeta}
\ee
Equation \eq{eq:zeta} establishes the important connection between the Dedekind $\eta$-function and
the respective series, that we substantially exploit below.

\subsection{Relation between the popcorn and Dedekind $\eta$ functions}

Consider an arbitrary quadratic form $Q'(m,n)$ with unit determinant. Since $d=1$, it can be
written in new parameters $\{a, b, c\} \to \{x = \frac{b}{c}, \eps = \frac{1}{c}\}$ as follows:
\be
Q'(m, n) = \frac{1}{\eps}(x m - n)^2 + \eps m^2
\label{q}
\ee
Applying the first Kronecker limit formula to the Epstein function with \eq{q} and $s = 1 + \tau$,
where $\tau \ll 1$, but finite one gets:
\be
\zeta(s,Q') = \frac{\pi}{s-1} + 2\pi \left(\gamma + \ln\sqrt{\frac{1}{4\eps}} - 2\ln|\eta(x+i
\eps)|\right) + O(s-1)
\label{zeta1}
\ee
On the other hand, one can make use of the $\eps$-continuation of the Kronecker $\delta$-function,
\eq{delt}, and assess $\zeta(1 + \tau, Q')$ for small $\tau \ll 1$ as follows:
\be
\zeta(1 + \tau, Q') \approx \frac{1}{\eps}\sum_{\{m,n\} \in \mathbb{Z}^2 \backslash \{0, 0\}}
\frac{\eps^2}{\left(x m - n\right)^{2} + \eps^2m^2} = \frac{2}{\eps}\lim_{N\to\infty}
\sum_{m=1}^{N} \sum_{n=1}^{N} \frac{1}{m^2}\delta\left(x - \frac{n}{m}\right) \equiv \theta(x);
\qquad x \in (0, 1)
\ee
where the factor 2 reflects the presence of two quadrants on the $\mathbb{Z}^2$-lattice that
contribute jointly to the sum at every rational points, while $\theta$ assigns 0 to all irrationals.
At rational points $\theta\left(\frac{p}{q}\right)$ can be calculated straightforwardly:
\be
\theta\left(\frac{p}{q}\right) = \frac{2}{\eps} \sum_{m|q}^{\infty} \frac{1}{m^2} =
\frac{\pi^2}{3 \eps q^2}
\label{eq:14}
\ee
Comparing \eq{eq:14} with the definition of the popcorn function, $g$, one ends up with the
following relation at the peaks:
\be
g\left(\frac{p}{q}\right) = \sqrt{\frac{3\eps}{\pi^2}\theta\left(\frac{p}{q}\right)}
\label{popc}
\ee
Eventually, collecting \eq{zeta1} and \eq{popc}, we may write down the regularization of the
popcorn function by the Dedekind $\eta (x + i\eps)|_{\eps \to 0}$ in the interval $0 < x < 1$:
\be
g(x) \approx \sqrt{-\frac{12\eps}{\pi} \ln|\eta(x+i\eps)| - o\left(\eps\ln\eps\right)} \Bigg|_{\eps
\to 0}
\label{result}
\ee
or
\be
-\ln|\eta(x+i\eps)|_{\eps \to 0} = \frac{\pi}{12\eps} g^2(x) + O(\ln\eps)
\label{result1}
\ee

\begin{figure}[ht]
\centerline{\includegraphics[width=13cm]{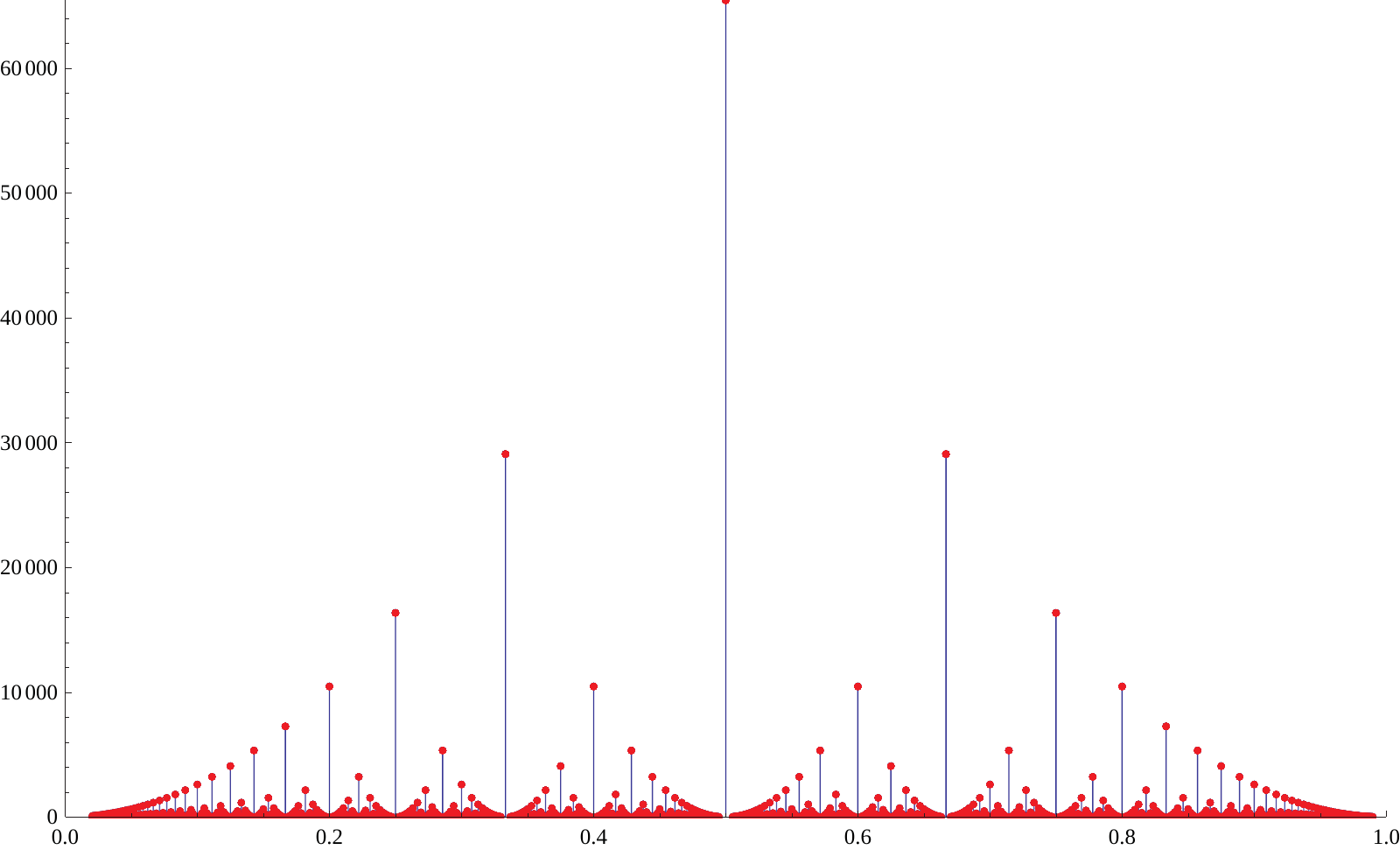}}
\caption{Plots of everywhere continuous $f_1(x) = -\ln|\eta(x+i\eps)|$ (blue) and discrete
$f_2(x) = \frac{\pi}{12\eps} g^2(x)$ (red) for $\eps = 10^{-6}$ at rational points in $0<x<1$.}
\label{fig:resfig}
\end{figure}

Note, that the asymptotic behavior of the Dedekind $\eta$-function can be independently derived
through the duality relation, \cite{avetisov}. However, such approach leaves in the dark the
underlying structural equivalence of the popcorn and $\eta$ functions and their series
representation on the lattice $\mathbb{Z}^2$. In the \fig{fig:resfig} we show two discrete plots of
the left and the right-hand sides of \eq{result1}.

Thus, the spectral density of ensemble of linear chains, \eq{eq:115}, in the regime $\eps \ll 1$
is expressed through the Dedekind $\eta$-function
as follows:
\be
\rho_\eps (\lambda) \approx \sqrt{-\frac{12\eps}{\pi} \ln\Bigg|\eta\left(\frac{1}{\pi}\arccos{\frac{\lambda}{2}}+i\eps\right)\Bigg|}
\ee

\subsection{Continued Fractions and Invariant Measure}

The behavior of the spectral density of our model is very similar to the one encountered in the
study of harmonic chains with binary distribution of random masses. This problem, which goes back
to F. Dyson \cite{dyson}, has been investigated by C. Domb \emph{et al.} in \cite{domb} and then
thoroughly discussed by T. Nieuwenhuizen and M. Luck in \cite{huisen}. Namely, consider a harmonic
chain of masses, which can take two values, $m$ and $M>m$:
\be
\left\{\begin{array}{cl}
m & \mbox{with the probability $1-f$} \medskip \\
M & \mbox{with the probability $f$}
\end{array}\right.
\label{co:1}
\ee
In the limit $M\to\infty$ the system breaks into "uniform harmonic islands", each consisting of $n$
light masses (in what follows we take $m=1$ for simplicity) surrounded by two infinitely heavy
masses. The probability of such an island is $(1-f)^2 f^n$. There is a straightforward mapping of
this system to the exponentially weighted ensemble of linear graphs discussed in the previous
section.

Many results concerning the spectral statistics of ensemble of sequence of random composition of
light and heavy masses were discussed in the works \cite{domb, huisen}. In particular, in
\cite{huisen} the integrated density of states, ${\cal N}(\lambda)$ has been written in the
following form:
\be
{\cal N}(\lambda) = \int_{-\infty}^{\lambda} \rho(\lambda')d\lambda' =
1-\frac{1-f}{f^2}\sum_{n=1}^{\infty} f^{\rm Int\left(\frac{n\pi}{\vartheta}\right)}
\label{co:20}
\ee
where the relation between $\lambda$ and $\vartheta$ is the following:
\be
\cos\vartheta=\frac{\sqrt{\lambda+1}}{2} \qquad \left(0<\vartheta<\frac{\pi}{2}\right)
\label{co:7}
\ee
Note that \eq{co:7} implies the asymmetry in the spectrum, $-1 < \lambda < 3$, that is related with
presence of additional "1" on the main diagonal of the adjacent matrix $B'$, obtained from the
matrix $B$ defined in \eq{eq:06}. Both matrices are deeply related to the spectrum of resonance
frequencies for the ensemble of linear chains with the distribution $Q_n$. Approximating the
Laplacian by $O(h^2)$-difference operator and rewriting the matrix eigenvalues problem as a
boundary problem for a grid eigenfunction $y_k(i)$, corresponding to $\lambda_{k, n}$, one has:
\be
\frac{y_k(i-1) - 2y_k(i) + y_k(i+1)}{h^2} + w^2_{k,n} y_k(i) = 0; \qquad 2 \le i \le n-1
\ee
Since the frequencies $w^2_{k, n}$ are independent of the choice of the specific matrix ($B$ or
$B'$), being the characteristics of the ensemble, one may write down the relation between the
eigenvalues $\lambda_{k,n}$ and $\lambda'_{k, n}$:
$$
\lambda_{k, n} = w^2_{k, n} + 2; \qquad \lambda'_{k, n} = w^2_{k, n} + 3
$$
Thus, the spectra of $B$ and $B'$ are linked by the shift $\lambda'_k = \lambda_k+1$ and spectrum
of $B'$ is asymmetric.

For $\lambda\to -1$ one gets ${\cal N}(\lambda)\to f/(1+f)$ which corresponds to the contribution
of the states $\lambda=-1$ at the bottom of the spectrum. At the upper edge of the spectrum (for
$\lambda\to 3^{-}$) one gets ${\cal N}(\lambda)\to 1$ which means that all states are counted.
Equation (\ref{co:20}) shows that the behavior near $\lambda=3$ (corresponding to $\vartheta=0$) is
dominated by the first term ($n=1$) of the series. One has:
\be
{\cal N}(\lambda) \simeq 1-\frac{1-f}{f^2}\exp\left(\frac{2\pi\ln f}{\sqrt{3-\lambda}}\right)
\label{co:21}
\ee
The expression (\ref{co:21}) signals the appearance of the Lifshitz singularity in the density of
states. A more precise analysis shows that the tail \eq{co:21} is modulated by a periodic function
\cite{huisen}.

Equation (\ref{co:20}) displays many interesting features. In particular, the function ${\cal
N}(\lambda)$ occurs in the mathematical literature as a generating function of the continued
fraction expansion of $\disp\frac{\pi}{\vartheta}$. Let us briefly sketch this connection following
\cite{borwein}. Consider the continued fraction expansion:
\be
\frac{\pi}{\vartheta}=\frac{1}{\disp c_0+ \frac{1}{\disp c_1+ \frac{1}{\disp c_2+\ldots}}}
\label{co:23}
\ee
where all $c_n$ are natural integers. Truncating this expansion at some level $n$, one gets a
coprime quotient, $\frac{p_n}{q_n}$ which converges to $\frac{\pi}{\vartheta}$ when $n\to \infty$.
A theorem of Borwein and Borwein \cite{borwein} states that the generating function, $G(z)$ of the
integer part of $\frac{\pi}{\vartheta}$
$$
G(z)=\sum_{n=1}^{\infty}z^{\rm Int\left(\frac{n\pi}{\vartheta}\right)}
$$
is given by the continued fraction expansion
$$
G(z)=\frac{z}{1-z} \frac{1}{\disp A_0+ \frac{1}{\disp A_1+ \frac{1}{\disp A_2+\ldots}}}
$$
where
$$
A_n(z)=\frac{z^{-q_n}- z^{-q_{n-2}}} {z^{-q_{n-1}}-1}
$$
and $q_n$ is the denominator of the fraction $\frac{p_n}{q_n}$ approximating the value
$\disp\frac{\pi}{\vartheta}$. In order to connect $G(z)$ with the integrated density, \eq{co:20},
it is sufficient to set $z=f$ and express ${\cal N}(\lambda)$ in terms of $G(z)$. Then, from the
equivalence discussed at the beginning of this section, the derivative of $N(\lambda)$ yields
the popcorn function.

The connection between the popcorn function to the invariant measure has been discussed in
different context in the work \cite{nechaev-comtet}.

\section{Conclusion}

We have discussed the number-theoretic properties of distributions appearing in physical systems
when an observable is a quotient of two independent exponentially weighted integers. The spectral
density of ensemble of linear polymer chains distributed with the law $f^L$ ($0<f<1$), where $L$ is
the chain length, serves as a particular example. In the sparse regime, namely at $f\to 1$, the
spectral density can be expressed through the discontinuous and non-differentiable at all rational
points, Thomae ("popcorn") function. We suggest a continuous approximation of the popcorn function,
based on the Dedekind $\eta$-function near the real axis.

Analysis of the spectrum at the edges reveals the Lifshitz tails, typical for the 1D Anderson
localization. The non-trivial feature, related to the asymptotic behavior of the shape of the
spectral density of the adjacency matrix, is as follows. The main, enveloping, sequence of peaks
$1-2-3-4-5...$ in the \fig{fig:02} has the asymptotic behavior $\rho(\lambda)\sim
q^{\pi/\sqrt{2-\lambda}}$ (at $\lambda\to 2^-$) typical for the 1D Anderson localization, however
any internal subsequence of peaks, like $2-6-7-...$, has the behavior $\rho'(\lambda)\sim
q^{\pi/|\lambda-\lambda_{cr}|}$ (at $\lambda\to \lambda_{cr}$) which is reminiscent of the Anderson
localization in 2D.

We would like to emphasize that the ultrametric structure of some collective observables, like the
spectral density of ensemble of sparse random Schr\"odinger-like operators, is ultimately connected
to number-theoretic properties of modular functions, like the Dedekind $\eta$-function, and
demonstrate the hidden $SL(2,Z)$ modular symmetry. We also pay attention to the connection of the
Dedekind $\eta$-function near the real axis to the invariant measures of some continued fractions
studied by Borwein and Borwein in 1993 \cite{borwein}.

The notion of ultrametricity deals with the concept of hierarchical organization of energy
landscapes \cite{mez,fra}. A complex system is assumed to have a large number of metastable states
corresponding to local minima in the potential energy landscape. With respect to the transition
rates, the minima are suggested to be clustered in hierarchically nested basins, i.e. larger basins
consist of smaller basins, each of these consists of even smaller ones, \emph{etc}. The basins of
local energy minima are separated by a hierarchically arranged set of barriers: large basins are
separated by high barriers, and smaller basins within each larger one are separated by lower
barriers. Ultrametric geometry fixes taxonomic (i.e. hierarchical) tree-like relationships between
elements and, speaking figuratively, is closer to Lobachevsky geometry, rather to the Euclidean
one.

\begin{acknowledgments}
We are very grateful to V. Avetisov, A. Gorsky, Y. Fyodorov and P. Krapivsky for many illuminating
discussions. The work is partially supported by the IRSES DIONICOS and RFBR 16-02-00252A grants.

\end{acknowledgments}

\begin{appendix}

\section{Asymptotics of $\eta(x+iy)$ at rational $x$ and $y\to 0$}

Here we shall obtain the asymptotics of $\eta(z)$ at $y\to 0^+$. Take into account the
relation of the Dedekind $\eta$ with Jacobi elliptic functions:
\be
\vartheta_1'(0,e^{\pi i z})=\eta^3(z)
\label{a:17}
\ee
where
\be
\vartheta_1'(0,e^{\pi i z})\equiv \frac{d\vartheta_1(u,e^{\pi i z})}{du}\bigg|_{u=0}= e^{\pi i
z/4}\sum_{j=0}^{\infty}(-1)^j (2j+1)e^{i\pi j(j+1)z}
\label{a:18}
\ee
Rewrite \eq{eq:15a} in terms of $\eta$-functions:
\be
\disp \left|\eta\left(\left\{\frac{m}{k}\right\} + iy\right)\right| = \frac{1}{\sqrt{k y}}
\left|\eta\left(\left\{\frac{n}{k}\right\}+\frac{i}{k^2y} \right)\right|
\label{a:19}
\ee
Applying \eq{a:17}--\eq{a:18} to \eq{a:19}, we get:
\be
\left|\eta\left(\left\{\frac{m}{k}\right\} + iy\right)\right|=\frac{1}{\sqrt{k y}}
\left|e^{\frac{\pi}{4} i \left(\left\{\frac{n}{k}\right\}+\frac{i}{k^2y}\right)}
\sum_{j=0}^{\infty}(-1)^j (2j+1)e^{i\pi j(j+1)\left(\left\{\frac{n}{k}\right\}+\frac{i}{k^2y}
\right)} \right|^{1/3}
\label{a:20}
\ee
Equation \eq{a:20} enables us to extract the leading asymptotics of $\left|\eta\left(\left\{
\frac{m}{k}\right\} + iy\right)\right|$ in the $y\to 0^+$ limit. Note that every term in the series
in \eq{a:20} for any $n\ge 1$ and for small $y$ (i.e. for $y\ll k^{-2}$) converges exponentially
fast, we have
\be
\left|\eta\left(\left\{\frac{m}{k}\right\} + iy\right)\right|_{y\to 0^+}\to\frac{1}{\sqrt{k y}}
e^{-\frac{\pi}{12 k^2y}};
\label{a:21}
\ee
Thus, we have finally:
\be
\sqrt{-\ln \left|\eta\left(\left\{\frac{m}{k}\right\} + iy\right)\right|}\Bigg|_{y\to 0^+}
\to\frac{1}{k}\sqrt{\frac{\pi}{12y}}
\label{a:23}
\ee
In \eq{a:23} we have dropped out the subleading terms at $y\to 0^+$, which are of order of $O(\ln
y)$.

\end{appendix}

\end{document}